\documentclass[12pt]{article}
\usepackage[T2A]{fontenc}
\usepackage{amsmath,amsthm}
\usepackage{amsfonts,amssymb}
\usepackage{mathrsfs}
\usepackage{graphicx}
\usepackage{hyperref}
\begin{document}

\begin{center}
{\bf Meta-intransitive Systems of Independent Random Variables and Fractals}

~

Alexey V. Lebedev\footnote{Faculty of Mechanics and Mathematics, Lomonosov Moscow State University, Moscow, Russia}
\end{center}

{\bf Abstract}

We construct  meta-intransitive systems of independent random variables of any finite
order from basic tuple of random variables which generalize intransitive dice.
Under this construction, the equality of some linear functional is preserved,
provided this equality hold for the basic variables.
This scheme is also extended to infinite-order systems, which are shown to be fractals, and for which
the similarity dimension can be evaluated. The problem of the upper bound for this dimension is posed.

~


~


~

~\\
Alexey V. Lebedev\\
avlebed@yandex.ru\\
ORCID 0000-0002-9258-0588

\newpage

\section{Introduction}

In theoretical studies and in practical applications, various superiority
relations between objects frequently feature the transitivity property: if $A$ is superior to~$B$ and $B$ is superior to~$C$,
then  $A$ is superior to~$C$.

However, other cases are also possible, For example,
in the classic finger-flashing game
rock--scissors--paper: rock beats scissors, scissors beats paper, paper beats rock, and not vice versa.
Thus, the relation between strategies in this game is intransitive.

Various aspects and numerous examples of intransitivity in superiority relations in nature, technology and society
have been addressed by Poddiakov, see \cite{PodVal, Pod2}, etc.

Among the recent works devoted to applications of intransitivity in biology,
we note the interesting studies on the life of ants in Puerto Rico \cite{Vander}
and soil bacterial networks \cite{Verdu}.

In what follows, we shall speak about intransitivity of the {\it stochastic precedence} relation
between random variables, namely, when one variable
is frequently smaller than another one.
In~\cite{Arc, Bol},  this relation was applied in
statistical analysis problems.
In~\cite{Shakhnov}, the reverse relation was applied in ranking problems
 (where it was called the dominance with respect to probability).

Probabilistically, the intransitivity problem was first considered by Trybula (jointly with Steinhaus),
who studied, as a~practical application, the strength of materials
\cite{Stein1959, Tryb1961}.
Assume that pairwise strength tests of iron bars of the same size and shape
are conducted. The bars are produced in three different factories (with different strength distribution in bars) and
that the bars are compared from the first and second, second and third, first and third factories.
In theory,  a~paradoxical situation may arise: bars from the first factory are ``worse''  (break more often)
than those from the second factory, bars from the second factory are ``worse'' than those from the third
factory, and bars from the third factory are ``worse'' than those from the first factory.

Recently \cite{Demler}, a~similar phenomenon was found to occur in medical statistics when comparing time to death and myocardial infarction in randomized controlled trial with three treatment arms.

The topic of intransitivity of stochastic precedence has gained popularity on the example of the so-called
\textit{intransitive} (\textit{nontransitive}) \textit{dice}.
Here one means sets of dice whose faces contain numbers  resulting in intransitive relations between the random variables.
It is assumed that two players choose dice from some set and throw them. The one wins who has a~larger number.
Intransitive dice  were popularized by Gardner\cite{Gard1,Gard2},  and have been extensively  studied.
Under this approach, the concept of a~die  is generalized\,---\,namely, one considers
dice with $n$~faces,  $n\ge 3$, labeled with arbitrary numbers $c_1,\dots, c_n$.
It is assumed that each face is equiprobable.
It can be also assumed without loss of generality that all the numbers are arranged in nondecreasing order:
$$
c_1\le c_2\le\dots\le c_n.
$$

There is extensive literature on this topic. As recent works, see, for example, \cite{Akin,Akin2,Boz,Con,Hazla,Kirk}.

The present paper continues author's studies on intransitivity began in~\cite{Leb2019_1,GorLeb2020,GorLeb2022,
PodLeb} (partially jointly with A.\,V.~Gorbunova  and A.\,N.~Poddiakov).
In these papers, the reader may find more details on the problems considered, historical account, and applications.

\section{Definitions and earlier results}

We first consider arbitrary sets of objects and relations defined on these objects.

A relation $\prec$ is an  {\it intransitive relation} if, for some objects $A$, $B$, $C$,
the relations $A\prec B$, $B\prec C$ do not imply that  $A\prec C$
(on the contrary, $C\prec A$ is possible).

A tuple of objects $A$, $B$, $C$ is called an {\it intransitive tuple} if on this tuple
the relation  $\prec$ is intransitive, i.e., if
$A\prec B\prec C\prec A$ (or in the reversed order).

A tuple of objects $A_1,\dots, A_m$, $m\ge 3$, will be said to be an intransitive
tuple of length~$m$ if
$$
A_1\prec A_2\prec\dots\prec A_m\prec A_1.
$$
So, an intransitive relation establishes a~cyclic order between objects
of an intransitive tuple. More involved generalizations are also possible.

The above relations and tuples of objects induce relations  and tuples of sets of objects.

Given sets   $M_1$, $M_2$, we say that $M_1\prec M_2$ if $A_1\prec A_2$ for all $A_1\in M_1$, $A_2\in M_2$.

Sets $M_1,\dots, M_m$, $m\ge 3$, form an intransitive tuple of length $m$ if
$$
M_1\prec M_2\prec\dots\prec M_m\prec M_1,
$$
which is equivalent to  saying that, for any choice of elements $A_1\in M_1,\dots, A_m\in M_m$,
these elements form an intransitive tuple (in the same order of traverse).

We will assume that  $X$ {\it stochastically precedes} $Y$ ($X\prec Y$) if
$${\bf E}{\rm\ sign}(Y-X)>0$$
or, what is the same,
$$
{\bf P}(X<Y)>{\bf P}(X>Y).
$$

We can also introduce the relation via the {\it preference function}
$$
\rho(Y,X)={\bf E}{\rm\ sign}(Y-X);
$$
in this case, $X\prec Y$ is equivalent to $\rho(Y,X)>0$. Here, the function
reflects not only the fact of
superiority  (or precedence), but also its measure.

If ${\bf P}(X=Y)=0$ (for example, random variables are independent and continuous,
or have disjoint ranges) then  $X\prec Y$ is equivalent to saying that
$$
{\bf P}(X<Y)>\frac{1}{2}.
$$
In this case, an intransitivity situation occurs if the random variables $X$, $Y$ and $Z$ satisfy
$$
P_{XYZ}= \min\{{\bf P}(X<Y),{\bf P}(Y<Z),{\bf P}(Z<X)\}>\frac{1}{2};
$$
here, not only the mere fact is of interest, but also the measure of intransitivity.

According to \cite{Tryb1961}, for independent random variables,
$$
\max_{X,Y,Z}P_{XYZ} = \frac{\sqrt{5}-1}{2}\approx 0.618,
$$
and the maximum is attained, for example, on the triplet
\begin{equation}\label{tripl0}
\begin{array}{l}
X=\left\{\begin{array}{ll}
1, & \mbox{ with probability } p,\\
4, & \mbox{ with probability } 1-p,
\end{array}\right.\\
Y=2,\\
Z=\left\{\begin{array}{ll}
0, & \mbox{ with probability } 1-p,\\
3, & \mbox{ with probability } p,
\end{array}\right.
\end{array}
\end{equation}
where
$$
p=\frac{\sqrt{5}-1}{2}.
$$
Hence
$$
{\bf P}(X<Y)={\bf P}(Y<Z)={\bf P}(Z<X)=\frac{\sqrt{5}-1}{2}.
$$

In \cite{Tryb1965,Usis}, combinations of $m\ge 3$ independent  random variables were studied.
Based of these studies, it was shown in~\cite{Bogd} that the  maximum of the probabilities
$$
P_{X_1\dots X_m}=\min\{{\bf P}(X_1<X_2),\dots,{\bf P}(X_{m-1}<X_m),{\bf P}(X_m<X_1)\}
$$
is as follows:
\begin{equation}\label{pmax}
\max\limits_{X_1,\dots, X_m}P_{X_1\dots X_m}=1-\left(4\cos^2\frac{\pi}{m+2}\right)^{-1},\quad m\ge 3.
\end{equation}
(later,  in~\cite{Komis} this result was reproved by a~geometric method).
This maximum is attained, for example, in the case of Efron dice (with $m=4$).

{\it Efron dice} (ED), which were invented by B.~Efron in the 1960s, is a~tuple of four
dice A, B, C, D with the numbers of the faces:
\begin{center}
    A: 0, 0, 4, 4, 4, 4;

    B: 3, 3, 3, 3, 3, 3;

    C: 2, 2, 2, 2, 6, 6;

    D: 1, 1, 1, 5, 5, 5.
\end{center}
In this case, the result of throwing each die from this set is grater than that for the succeeding
(in the circle) die with probability exceeding $1/2$:
$$
{\bf P}(A>B)={\bf P}(B>C)={\bf P}(C>D)={\bf P}(D>A)=\frac23,
$$
To get the reciprocal relations, one may consider dice in the reversed order (D, C, B, A).

For  $m\ne 4$, the maximum in \eqref{pmax} is not attained on dice, since this maximum is
an irrational number.

\smallskip

Let us give some other examples of intransitive dice.

{\it Classical intransitive 3-sided dice} (CID) are obtained from ancient Chinese Lo Shu magic square:

\begin{center}
\begin{tabular}{|c|c|c|}
\hline
4 & 9 & 2\\
\hline
3 & 5 & 7\\
\hline
8 & 1 & 6\\
\hline
\end{tabular}
\end{center}
Namely, by ordering the numbers in the  rows in an increasing order, we set
\begin{center}
A: 2, 4, 9;\\
B: 3, 5, 7;\\
C: 1, 6, 8.\\
\end{center}
Hence
$$
{\bf P}(A<B)={\bf P}(B<C)={\bf P}(C<A)=5/9.
$$

Moon and Moser \cite{MM} seem to be the first to observe that the magic square
gives an intransitive tuple for triples of numbers (columns were picked from the square in \cite{MM},
and rows, in~\cite{Gard2}). We will follow the second variant.
Note that \cite{MM,Gard2} were not concerned with a~game of dice, but rather with
competitions of sports teams, but this is not important.
Here, each sum of numbers on faces is~15, and the mean is~5.

There also exists a~modification of Efron dice with equal means.
For more involved examples, see, for example, \cite{Akin,Akin2,Boz,Con}, etc.
Examples of dice for which not only the means are equal, but the variances are also equal,
were recently presented in~\cite{Yak} for $n\ge 8$.

{\it Simplest intransitive 3-sided dice} (SID) are obtained from triplet (\ref{tripl0}) with $p=2/3$:
\begin{center}
A: 1, 1, 4;\\
B: 2, 2, 2;\\
C: 0, 3, 3.\\
\end{center}
We have
$$
{\bf P}(A<B)=2/3,\quad {\bf P}(B<C)=2/3,\quad {\bf P}(C<A)=5/9.
$$
In this case, the first player selects die~C, the second player selects die~A, and
the guaranteed probability of winning for the second player is $P_{ABC}=5/9$.
The sum of the numbers on the faces is~ 6, the mean is~2.

The meta-intransitivity (of type I) was introduced in \cite{PodLeb}. Let us clarify this concept.

{\it A meta-intransitive system} is a~set equipped with an intransitive relation of elements with the
following properties.

1. The set is partitioned into groups  (subsets), which may have different nesting level.

2. The set is partitioned into first-level groups, which form an intransitive tuple,
which, in turn, may either contain intransitive sets of elements or
be partitioned into second-level groups which form intransitive tuples, and so~on.

The {\it order} of a meta-intransitive is defined as the maximal nesting level of  groups (possibly infinite).

As an illustrative example, we mention the nested Borromean rings \cite{Bor}.
For some examples, see also \cite{PodLeb,Pod2023}. Meta-intransitive system of second order
presented in Fig.~1 (schematically).

\begin{figure}[htb]
\centering
\includegraphics[height=0.3\textheight]{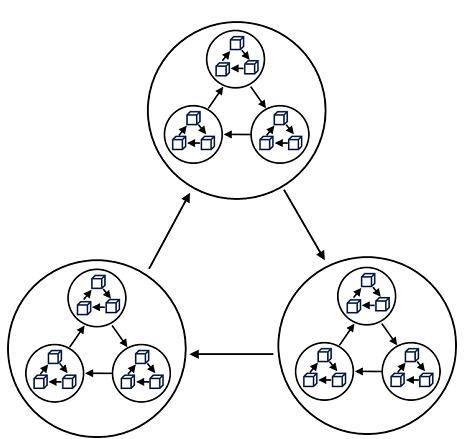}
\caption{Meta-intransitive system of second order \cite{Pod2023}.}
\end{figure}

Our next purpose is to construct and study meta-intransitive systems of random variables and fractals
based on a~generalization of intransitive dice.

\section{Main notation and assumptions}

We note, first of all, that with any random variable~$X$
one can associate the inverse distribution function, also known as the
{\it quantile function}:
$$
x(u)=\inf\{ x: {\bf P}(X\le x)\ge u\},\quad u\in (0,1),
$$
and vice verse, from this function one can construct a~random variable with the same distribution, by taking
$x(U)$, where the random variable $U$ is uniformly distributed  on $[0,1]$.

The class of quantile functions ${\cal Q}$ is the class of all c\`ad functions which nondecrease on  $(0,1)$.
It is clear that this class is closed  with respect to addition of functions and multiplication by
positive numbers. Let $M\subset {\cal Q}$. For each  $a\in {\cal Q}$ and any number $b>0$, consider the set
\begin{equation}\label{abM}
a+bM=\{a+bx: x\in M\}\subset {\cal Q}.
\end{equation}

Note that in the case of an $n$-sided die with numbers  $c_1\le c_2\le\dots\le c_n$ on faces, we get
\begin{equation}\label{xdice}
x(u)=c_{[nu]+1},
\end{equation}
i.e., we get the value $c_1$ on  $[0,1/n)$, the value $c_2$ on~$[1/n,2/n)$, and so~on.

Let us define the preference function on $x_1,x_2\in {\cal Q}$ by
\begin{equation}\label{rho12}
\rho_q(x_1,x_2)=\int_0^1\int_0^1 {\rm\ sign}(x_1(u_1)-x_2(u_2))\,du_1\,du_2.
\end{equation}
Hence, for the independent random variables $X_1$ and $X_2$ with such quantile functions, we have
\begin{equation}\label{rhoeq}
\rho(X_1,X_2)=\rho_q(x_1,x_2).
\end{equation}

For the piecewise continuous function $f:[0,1]\to {\mathbb R}$, consider the linear functional
$$
J_f(X)=\int_0^1 f(u)x(u)\,du
$$
(assuming that the integral converges). In the simplest case $f\equiv 1$, we have $J_1(X)={\bf E}X$.
Such functionals appear in the actuarial theory  (Wang's premium principle)
\cite{Wang}, financial mathematics (weighted VaR) \cite{Chern},
modern extreme value theory~\cite{Ivanov}, etc.

In the case of $n$-sided dice with numbers on faces $c_1\le c_2\le\dots\le c_n$, we get
\begin{equation}\label{Jdice}
J_f(X)=\sum_{i=1}^n\left(\int_{\frac{i-1}{n}}^{\frac{i}{n}}f(u)\,du\right)x\left(\frac{i-1}{n}\right)=
\sum_{i=1}^n\left(\int_{\frac{i-1}{n}}^{\frac{i}{n}}f(u)\,du\right)c_i.
\end{equation}

Consider a~basic tuple of independent random variables $X_1^0,\dots,X_m^0$, $m\ge 3$, such that
\begin{equation}\label{X0}
X_1^0\prec X_2^0\prec\dots \prec X_m^0\prec X_1^0
\end{equation}
and
\begin{equation}\label{Rr}
\inf_{u_{1,2}\in (0,1)}\min_{i_1\ne i_2}|x_{i_1}^0(u_1)-x_{i_2}^0(u_2)|=r>0,\quad
\sup_{u_{1,2}\in (0,1)}\max_{i_1\ne i_2}|x_{i_1}^0(u_1)-x_{i_2}^0(u_2)|=R<\infty.
\end{equation}
From condition \eqref{Rr} it follows that the random variables are bounded and have disjoint ranges (as in the above examples of
dice). They characteristics are given in the table.

\begin{center}
\begin{tabular}{|l|c|c|c|c|}
\hline
Dice & $m$ & $n$ & $r$ & $R$ \\
\hline
ED & 4 & 6 & 1 & 6 \\
\hline
CID & 3 & 3 & 1 & 8 \\
\hline
SID & 3 & 3 & 1 & 4 \\
\hline
\end{tabular}
\end{center}

However, in the general case,   $X_1^0,\dots,X_m^0$
may fail to be discrete and assume only finite number of values.
For example, integer-valued dice can be augmented with continuous random variables on $[0,\Delta]$, where $0<\Delta<1$.

We set
\begin{equation}\label{usl1}
\lambda\ge 1+\frac{R}{r},\quad \varepsilon=\frac{1}{\lambda}.
\end{equation}

We let $A^{(k)}$, $k\ge 0$, denote the sets of quantile functions, which are called
{\it generations} (of functions). These functions have $k$-dimensional indexes.

Consider the following algorithm for construction of the generations.

1. $A^{(0)}$ contains one function $x\equiv 0$ with empty index.

2. Each quantile function $x\in A^{(k-1)}$ generates the quantile functions
$$
x(u)+\varepsilon^k x_1^0(u),\dots,x(u)+\varepsilon^k x_m^0(u),
$$
which are given indexes by affixing the numbers $1,\dots,m$, respectively, to the right of the mother number index.
This gives us the set $A^{(k)}$.

Now $A^{(k)}$ contains $m^k$ functions with  indexes $i_1i_2\dots i_k\in\{1,\dots, m\}^k$.

We let $A^{(k)}_{i_1\dots i_p}$, $1\le p\le k-1$, $k\ge 2$, denote the subsets
$A^{(k)}$ of functions such that their indexes start with  $i_1\dots i_p$ (common ancestors).
Hence $A^{(k)}$, $k\ge 2$, is partitioned into $A^{(k)}_{i_1}$; the subsets  $A^{(k)}_{i_1}$,
$k\ge 3$, are partitioned into  $A^{(k)}_{i_1i_2}$, and so on. So, we have nesting up to the level $k-1$.

Consider the tuple of random variables $U_{i_1\dots i_k}$, $i_1i_2\dots i_k\in\{1,\dots, m\}^k$,
which are independent and uniformly distributed  on $[0,1]$. We also set
$$
X_{i_1\dots i_k}=x_{i_1\dots i_k}(U_{i_1\dots i_k}).
$$

From these random variables, we compose the set $B^{(k)}$ and the subsets
$B^{(k)}_{i_1\dots i_p}$, $1\le p\le k-1$, $k\ge 2$, of random variables, whose indexes
start with $i_1\dots i_p$. For these random variables, we have the same partitioning and nesting
as for the functions in~$A^{(k)}$.
In view of (\ref{rhoeq}), they satisfy the same relations~$\prec$.

It remains to verify the intransitivity of tuples and subsets  of functions in~$A^{(k)}$.
This will imply that the tuples and subsets of random variables in~$B^{(k)}$ are intransitive.

\section{Main results for generations}

By construction, each quantile function from $A^{(k)}$, $k\ge 1$,
with index $i_1i_2\dots i_k\in \{1,\dots, m\}^k$ is described by the formulas
\begin{equation}\label{xik}
x_{i_1i_2\dots i_k}(u)=\sum_{l=1}^k\varepsilon^l x^0_{i_l}(u)=
\sum_{j=1}^m\left(\sum_{l=1}^k \varepsilon^l \delta_{i_lj}\right)x^0_j(u),
\end{equation}
where $\delta_{ij}$ is the Kronecker delta.

Note that in view of (\ref{xdice}) from the basic tuple $n$-sided dice
one gets only $n$-sided dice (in the sense of constancy of quantile functions
on the intervals $[(i-1)/n,i/n)$, $1\le i\le n$).

{\bf Theorem 1.} {\it For all functions $x',x''\in A^{(k)}$, $k\ge 1$, with different indexes
$i'_1i'_2\dots i'_k$ and $i''_1i''_2\dots i''_k$, we have  $\rho_q(x',x'')=\rho_q(x^0_{i'_\nu},x^0_{i''_\nu})$,
where $\nu=\min\{l:i'_l\ne i''_l\}$.
}

\begin{proof}\rm For all $u',u''\in (0,1)$, consider the difference
$$
\begin{array}{c} \displaystyle
x'(u')-x''(u'')=\sum_{l=1}^k \varepsilon^l\left(x_{i'_l}^0(u')-x_{i''_l}^0(u')\right)=\\ \displaystyle
=\varepsilon^\nu \left(\left(x^0_{i'_\nu}(u')-x^0_{i''_\nu}(u'')\right)+\sum_{l=1}^{k-\nu}\varepsilon^l
\left(x_{i'_{\nu+l}}^0(u')-x_{i''_{\nu+l}}^0(u')\right)\right).
\end{array}
$$
By \eqref{Rr}, we have
$$
\left|x^0_{i'_\nu}(u')-x^0_{i''_\nu}(u'')\right|\ge r
$$
and
$$
\sum_{l=1}^{k-\nu}\varepsilon^l \left(x_{i'_{\nu+l}}^0(u')-x_{i''_{\nu+l}}^0(u')\right)\le
R\varepsilon\frac{1-\varepsilon^{k-\nu}}{1-\varepsilon}<R\frac{\varepsilon}{1-\varepsilon}=\frac{R}{\lambda-1}\le r,
$$
by (\ref{usl1}). Hence
$$
\mathrm{sign} (x'(u')-x''(u''))=\mathrm{sign} (x^0_{i'_\nu}(u')-x^0_{i''_\nu}(u'')),
$$
which proves the theorem.
\end{proof}

\smallskip

The meaning of the theorem is as follows. in each generation, the functions
preserve the relations between their nearest noncommon ancestors.
If functions are in some relation, then all their descendants are in the same relation.
And these descendants form kindred groups (in subsequent generations).

{\bf Corollary 1.} {\it There is a one-one correspondence between the indexes and the quantile functions
of the same generation.}

\begin{proof}\rm
It is clear that to each index there corresponds a~quantile function by~(\ref{xik}).
Assume that there exist two different indexes $i'_1i'_2\dots i'_k$ and $i''_1i''_2\dots i''_k$
that generate the functions~$x'$ and~$x''$ via~\eqref{xik},
and which are equal: $x'=x''$. Hence  $\rho_q(x',x'')=0$.
But by Theorem~1 we have $\rho_q(x',x'')=\rho_q(x^0_{i'_\nu},x^0_{i''_\nu})\ne 0$,
which is a~contradiction. Hence to each quantile function of generation
there corresponds precisely one index.
\end{proof}

{\bf Corollary 2.} {\it The generation of quantile functions $A^{(k)}$, $k\ge 2$,
is
a~meta-intransitive system of order $k-1$, where as groups
of the $p$th nesting level one may take  $A^{(k)}_{i_1\dots i_p}$, $1\le p\le k-1$, $k\ge 2$.}

\begin{proof}\rm  The sets  have a~nesting structure by construction and Corollary~1. Let us
verify the intransitivity.

The set $A^{(k)}_{i_1\dots i_p}$ is partitioned into the sets $A^{(k)}_{i_1\dots i_{p-1}j}$, $1\le j\le m$.
We choose functions $x^{(j)}\in A^{(k)}_{i_1\dots i_{p-1}j}$. Hence by Theorem~1 we have
$$
\rho_q(x^{(j_1)},x^{(j_2)})=\rho_q(x^0_{j_1},x^0_{j_2}),\quad j_1\ne j_2.
$$
Hence $x^{(j)}$ are in the same relations as the functions from the basic tuple $x^0_j$,
i.e., they form an intransitive tuple of length~$m$.
Therefore, by definition, the sets $A^{(k)}_{i_1\dots i_{p-1}j}$, $1\le p\le k-1$, form an intransitive
tuple of length~$m$ for sets.

The set $A^{(k)}_{i_1\dots i_{k-1}}$ consists of the  functions $x_{i_1\dots i_{k-1}j}$, $1\le j\le m$,
for which, by Theorem~1, we have
$$
\rho_q(x_{i_1\dots i_{k-1}j_1},x_{i_1\dots i_{k-1}j_2})=\rho_q(x^0_{j_1},x^0_{j_2}),\quad j_1\ne j_2.
$$
Hence these functions are in the same relations as those from the basic tuple $x^0_j$,
i.e., they form an intransitive tuple of length $m$.
\end{proof}

{\bf Corollary 3.} {\it The generation of the random variables $B^{(k)}$, $k\ge 2$, is
a~meta-intransitive  system of order~$k-1$, where as groups of the $p$th nesting level
one may take $B^{(k)}_{i_1\dots i_p}$, $1\le p\le k-1$, $k\ge 2$.}

This result follows  from Corollary~2 by (\ref{rhoeq}).

By  $\sqcup$ we denote the disjoint union of sets. Form  (\ref{abM})
and (\ref{xik}) we have the recurrence formula:
\begin{equation}\label{rec1}
A^{(k)}=\bigsqcup_{i=1}^m A^{(k)}_i,\quad A^{(k)}_i=\varepsilon x^0_i+\varepsilon A^{(k-1)},\quad 1\le i\le m,
\quad k\ge 1.
\end{equation}

Below, this formula will be generalized to the case $k\to\infty$.

{\bf Proposition 1.} {\it If $J_f(X^0_i)=K$ for all $1\le i\le m$, then,  for all $X\in B^{(k)}$, $k\ge 1$,
\begin{equation}\label{JC}
J_f(X)=K\varepsilon\frac{1-\varepsilon^k}{1-\varepsilon}.
\end{equation}
}

\begin{proof}\rm
Let $X$ have an arbitrary index $i_1i_2\dots i_k$.
Since the functional is linear with respect to the function~$x$,
it follows from representation \eqref{xik}  that
$$
J_f(X)=\sum_{l=1}^k \varepsilon^l J_f(X^0_{i_l})=K\sum_{l=1}^k\varepsilon^l=K\varepsilon\frac{1-\varepsilon^k}{1-\varepsilon}.
$$
\end{proof}

For example, if the means of all dice in the basic tuple are equal  (as in the case of CID and SID),
their the means of all dice
of the same generation will also be equal.

As noted above, from $n$-sided dice one can obtain only $n$-sided dice.
With 3-sided
dice we associate the points $(x(0),x(1/3),x(2/3))$ of the Euclidean space. Then the equation
$J_f(X)=K$ by (\ref{Jdice}) defines a~plane. If all dice of the basic tuple (as space points)
lie in some plane, then all
dice of the same generation  also lie in some plane (which a~translation of the original one).
In the case of a~basic tuple  of three
3-sided dice, all the dice
of any generation always lie in some plane, because there is always a~plane passing through any
three points.

{\bf Example.} Consider an SID as a basic tuple. By (\ref{usl1}),\enskip  $\lambda\ge 5$.
Setting $\lambda=5$, we have $\varepsilon=0.2$. We construct  $A^{(3)}$, and associate with each function
the point $(x(0),x(1/3),x(2/3))$. The resulting set ${\tilde A}^{(3)}$ consists of  27 points (see Fig.~2).

\begin{figure}[htb]
\centering
\includegraphics[height=0.45\textheight]{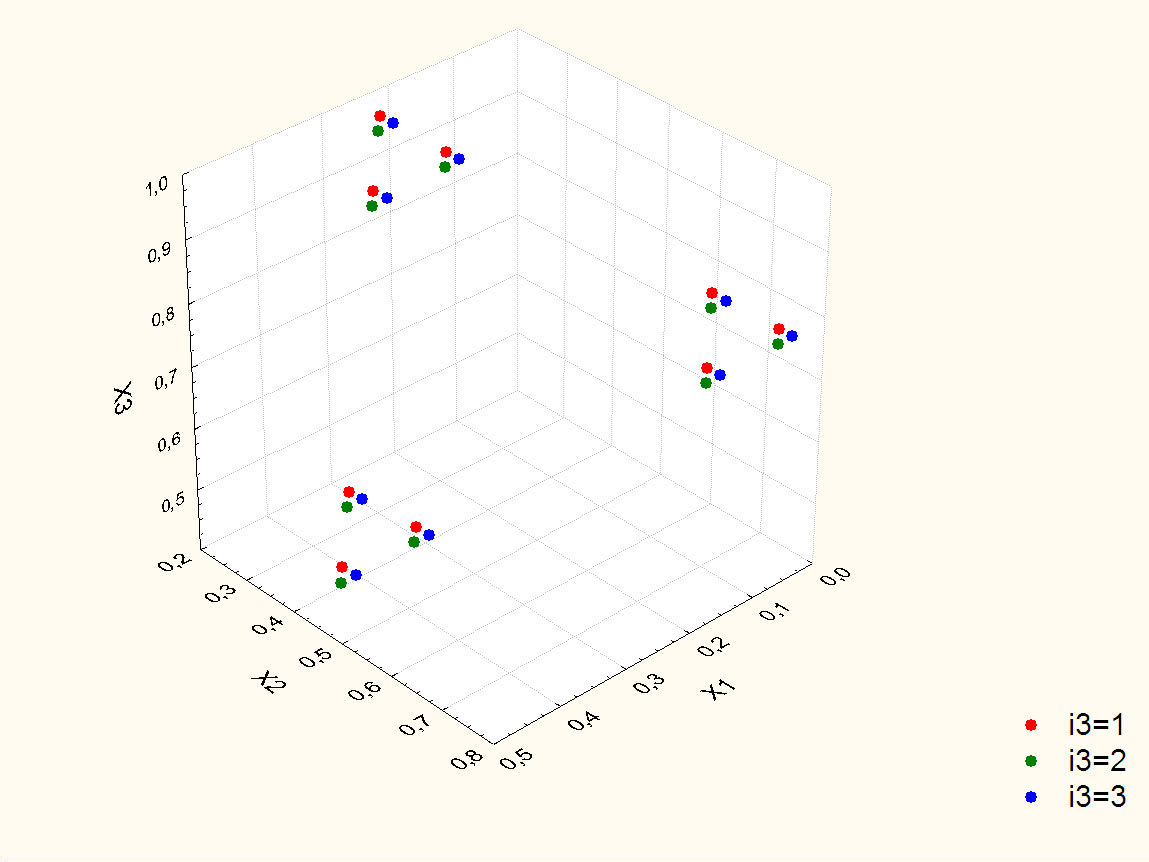}
\caption{The set  ${\tilde A}^{(3)}$ for SID (27 points), $\lambda=5$.}
\end{figure}

The figure is rotated for better clarity.
It is seen that all points lie on the same plane.
In the projection to the figure plane, the dice ``beat'' (dominate) each other in the clockwise
direction relative to the  center.
In the figure, the fractal nature of the set is clearly visible, which suggests the idea
of making $k\to\infty$ in our construction.

\section{Fractals}

Consider the following conditions, which is stronger than (\ref{usl1}):
\begin{equation}\label{usl2}
\lambda>1+\frac{R}{r},\quad \varepsilon=\frac{1}{\lambda}.
\end{equation}

We define the {\it infinite generation}
$A^{(\infty)}$ as the set of functions with infinite indexes $i_1\dots i_l\dots$ defined by the series
\begin{equation}\label{xinfty}
x_{i_1\dots i_l\dots}(u)=\sum_{l=1}^\infty\varepsilon^lx^0_{i_l}(u)=
\sum_{j=1}^m\left(\sum_{l=1}^\infty \varepsilon^l \delta_{i_lj}\right)x^0_j(u),
\end{equation}
where $\delta_{ij}$ is the Kronecker delta.

This is a denumerable set.

We let $A^{(\infty)}_{i_1\dots i_p}$, $p\ge 1$, denote the subsets
$A^{(\infty)}$ of dice such that their indexes start from $i_1\dots i_p$.
For these subsets, we have the same splitting and nesting
as in the finite case, but without constraints on~$p$.

We also introduce the tuple of
random variables $U_{i_1\dots i_l\dots}$, $1\le i_l\le m$,
jointly independent and uniformly distributed  on $[0,1]$. We also set
$$
X_{i_1\dots i_l\dots}=x_{i_1\dots i_l\dots}(U_{i_1\dots i_l\dots}).
$$

From the random variables set $B^{(\infty)}$ and subsets
$B^{(\infty)}_{i_1\dots i_p}$ we compose the random variables show indexes
start with  $i_1\dots i_p$. For these random variables, we have the same partitioning and nesting
as for the functions in~$A^{(\infty)}$. In view of (\ref{rhoeq}), they satisfy the same relations~$\prec$.

{\bf Theorem 2.}
{\it For all two functions from $x',x''\in A^{(\infty)}$ with different indexes
$i'_1\dots i'_l\dots$ and $i''_1\dots i''_l\dots$, we have  $\rho(x',x'')=\rho(c_{i'_\nu},c_{i''_\nu})$,
where $\nu=\min\{l:i'_l\ne i''_l\}$.
}

\begin{proof}\rm
For all $u',u''\in (0,1)$, consider the difference
$$
\begin{array}{c} \displaystyle
x'(u')-x''(u'')=\sum_{l=1}^\infty \varepsilon^l\left(x_{i'_l}^0(u')-x_{i''_l}^0(u')\right)=\\ \displaystyle
=\varepsilon^\nu \left(\left(x^0_{i'_\nu}(u')-x^0_{i''_\nu}(u'')\right)+\sum_{l=1}^\infty\varepsilon^l
\left(x_{i'_{\nu+l}}^0(u')-x_{i''_{\nu+l}}^0(u')\right)\right);
\end{array}
$$
Here, using  \eqref{Rr}, we have
$$\left|x^0_{i'_\nu}(u')-x^0_{i''_\nu}(u'')\right|\ge r$$
and
$$
\sum_{l=1}^\infty\varepsilon^l \left(x_{i'_{\nu+l}}^0(u')-x_{i''_{\nu+l}}^0(u')\right)\le
R\frac{\varepsilon}{1-\varepsilon}=\frac{R}{\lambda-1}<r,
$$
by(\ref{usl2}). So,
$$
\mathrm{sign} (x'(u')-x''(u''))=\mathrm{sign} (x^0_{i'_\nu}(u')-x^0_{i''_\nu}(u'')),
$$
which proves the theorem.
\end{proof}

{\bf Corollary 4.} {\it There is a one-one correspondence  between the indexes and the quantile functions $A^{(\infty)}$.
}

{\bf Corollary 5.} {\it The infinite generation of quantile functions $A^{(\infty}$ is a system of infinite order, where
as groups of $p$th nesting level
one may take  $A^{(k)}_{i_1\dots i_p}$, $p\ge 1$.

}
{\bf Corollary 6.} {\it The infinite generation of random variables $B^{(\infty)}$
is a~meta-intransitive system of infinite order, where as groups of the $p$th nesting level
one may take $B^{(k)}_{i_1\dots i_p}$, $p\ge 1$.
}

Corollaries 4--6 are proved similarly to Corollaries 1--3.

{\bf Proposition 2.} {\it If $J_f(X^0_i)=K$ for all $1\le i\le m$, then  for all $X\in B^{(\infty)}$,
\begin{equation}\label{JC2}
J_f(X)=K\frac{\varepsilon}{1-\varepsilon}.
\end{equation}

}

The proof of Proposition 2 is similar to that of Proposition~1.

We also have the recurrence formula
\begin{equation}\label{rec2}
A^{(\infty)}=\bigsqcup_{i=1}^m A^{(\infty)}_i,\quad A^{(\infty)}_i=\varepsilon x^0_i+\varepsilon A^{\infty},
\quad 1\le i\le m,
\end{equation}
which is similar to~\eqref{rec1}.
So, the set $A^{(\infty)}$ consists of  $m$~parts similar to the entire set with similarity coefficient
$\varepsilon=1/\lambda$. Hence by \cite{Mand}, this is a~{\it fractal} with {\it similarity dimension}
$$
d=\frac{\ln m}{\ln\lambda}.
$$

By condition \eqref{usl2}, we get the boundary for dimension with known~$m$, $r$, and~$R$:
$$
d<d_{\sup}=\frac{\ln m}{\ln (1+R/r)}.
$$
This calls for the problem of the behavior of $d_{\sup}$.

Considering the characteristics of the above problems, we have the table.

\begin{center}
\begin{tabular}{|l|l|}
\hline
Dice &  $d_{\sup}$ \\
\hline
ED & $(\ln 4)/(\ln 7)\approx 0.7124$ \\
\hline
SID & $(\ln 3)/(\ln 5)\approx 0.6826$ \\
\hline
CID & $(\ln 3)/(\ln 9)=0.5$ \\
\hline
\end{tabular}
\end{center}

The maximal value is attained for the Efron dice. In all cases $d<1$, and hence
by \cite{Mand} the fractals are subsumed into the  {\it fractal dust} category.

This suggests the interesting question about the boundary for the dimension in the general case.
This question is still open.

\end{document}